\documentclass{article}
\usepackage{multirow}
\usepackage[utf8]{inputenc}
\usepackage[T1]{fontenc}
\usepackage{amsmath,amsfonts,amssymb,amsthm,amscd,tipa,color,graphicx,tikz-cd,longtable}
\usepackage[english]{babel}
\usepackage[a4paper,margin=0.75in,top=4.5cm,bottom=3.5cm]{geometry}
\usepackage[all]{xy}
\usepackage{multirow}
\usepackage{float}

\newcommand{\Q}{\mathbb{Q}}

\newcommand{\C}{\mathbb{C}}

  \DeclareMathOperator{\pgl}{PGL}

\theoremstyle{plain}
\begingroup

\newtheorem{theo}{Theorem}[section]
\newtheorem*{theo*}{Theorem}

\begin{document}
\title{Representations of Deligne-Mostow lattices into $\pgl(3,\C)$. Part II}

\author{E. Falbel, I. Pasquinelli and A. Ucan-Puc}

\maketitle

\abstract{We complete the classification of type preserving representations of  Deligne-Mostow lattices with 3-fold symmetry into $\pgl(3,\C)$ started in \cite{ELA}. In particular, we show local rigidity for all the representations where the generators we chose are of the same type as the generators of the Deligne-Mostow lattices. We use formal computations in 
 SAGE and MAPLE to obtain the results.
The code files are available on GitHub (\cite{GH}). 

\section{Introduction}

This is a sequel to \cite{ELA}, where we started the classification of representations (up to conjugation) of Deligne-Mostow lattices into $\pgl(3,\C)$.

The general goal of this project is to describe the space of all representations of a Deligne-Mostow group $\Gamma$, into $\pgl(3,\C)$ up to conjugation.  
In other words, we want to study the space
$Hom(\Gamma,\pgl(3,\C))/\pgl(3,\C)$.  
In \cite{ELA} we studied Deligne-Mostow groups with 3-fold symmetry of type one (see section \ref{section:DM} for the definition of those groups).  The main result is a local rigidity theorem for representations preserving types of our choice of generators.  
The generators are elliptic elements and here, preserving type means that we are selecting those representations for which the image of the generators are either regular elliptic or complex reflections as they are in the Deligne-Mostow.

In this sequel we study the remaining Deligne-Mostow groups with 3-fold 
symmetry which are classified as types two, three and four according to the form of their presentations (see again section \ref{section:DM} for definitions).

The result is parallel to the type one case.  Namely, we prove local rigidity of representations with the condition that generators of our particular presentation 
have images of the same type as the Deligne-Mostow group.  Furthermore, we classify all these representations and give the expressions of the generators in their appropriate extension fields which are all cyclotomic. The tables in this paper contain, for each group, the number of representations, the number of Galois orbits and the field extension where the group elements can be defined.   We also detail which representations are irreducible and list the non-compact ones which factor through the orbifold fundamental group of a compactification. 
The detailed list of all representations
are given in code files with  computations with Groebner basis packages in SAGE and MAPLE  which are available on GitHub (\cite{GH}).  

The second author acknowledges the support of the FSMP and of the EPSRC grant EP/T015926/1. 
The third author acknowledges the support of the COnsejo NAcional de Ciencia Y Tecnolog\'ia, project 769037.
 
\section{Deligne-Mostow lattices with three-fold symmetry}
\label{section:DM}
In this section we recall the definition of the Deligne-Mostow lattices in $PU(2,1)$ (\cite{mostow},\cite{delignemostow}, see also \cite{ELA}). 
Let a \emph{ball 5-tuple} $\mu=(\mu_1, \dots, \mu_5)$ be a set of $5$ real numbers in $(0,1)$ whose sum is 2.
Deligne and Mostow define a  \emph{hypergeometric function}  using these parameters,
which gives rise to  the \emph{monodromy action} which defines  lattices in $PU(2,1)$.

The Deligne-Mostow lattices in $PU(2,1)$ with \emph{3-fold symmetry} correspond  
 to cases where  $3$ of the $5$ elements of the ball $5$-tuple are equal. These lattices depend on 2 parameters and 
a usual choice for the parameters is the orders of some of the generators. 
We  identify the lattices using a pair $(p,k)$ which corresponds to the ball 5-tuple $(1/2-1/p, 1/2-1/p, 1/2-1/p, 1/2+1/p-1/k, 2/p+1/k)$. 
A fundamental domain and a presentation for each of these lattices can be found in \cite{Irene}.  There are four types of lattices with 3-fold symmetry.  The types correspond to different complexities of the combinatorics of a fundamental domain.

%In this work we will only look at some of the lattices with 3-fold symmetry. 
%For all of those, one can find a construction for a fundamental domain and a presentation in \cite{Irene}.
%In principle the 2-dimensional Deligne-Mostow lattices depend on the 5 elements of the ball 5-tuple. 
%The condition on the sum of the elements of the ball 5-tuple reduces the parameters to 4. 
%Of these 4, 3 are equal, so the lattices only depend on 2 parameters. 
%One common choice for the parameters is to choose the orders of some of the generators. 
%We will hence identify the lattices using a pair $(p,k)$ which corresponds to the ball 5-tuple $(1/2-1/p, 1/2-1/p, 1/2-1/p, 1/2+1/p-1/k, 2/p+1/k)$. 

%The Deligne-Mostow lattices in $PU(2,1)$ with 3-fold symmetry are divided in 4 \emph{types} according to the ranges of $p$ and $k$. 
%The type determines the presentation, the volume formula and the combinatorics of a fundamental domain. 

%In terms of $p$ and $k$ they are characterised by 
%\begin{align*}
%0<p &\leq 6, & k \leq \frac{2p}{p-2}.
%\end{align*}

We will now describe the presentations we use all 3-fold symmetry lattices.

The lattices of type one were treated in \cite{ELA}.
We include their presentation here for completeness.
A 3-fold Deligne-Mostow lattice of type one has the following presentation:

\begin{equation*}
\left\langle J, P, R_1, R_2: J^3=R_1^p=R_2^p=(P^{-1}J)^k=Id, R_2=PR_1P^{-1} =J R_1 J^{-1}, P=R_1 R_2\right\rangle,
\end{equation*}

By taking $R_2=J R_1 J^{-1}$ and $P=R_1 R_2,$ it is possible to reduce the presentation to

\begin{equation}
\left\langle J, R_1: J^3=R_1^p=(R_1J)^{2k}=R_1 J R_1J^2 R_1 J R_1^{p-1} J^2 R_1^{p-1} J R_1^{p-1} J^2\right\rangle.
\end{equation}

A 3-fold Deligne-Mostow lattice of type two has the following presentation:

\begin{equation*}
\left\langle J, P, R_1, R_2: J^3=R_1^p=R_2^p=(P^{-1}J)^k=(R_2 R_1 J)^l=Id, R_2=PR_1P^{-1} =J R_1 J^{-1}, P=R_1 R_2\right\rangle,
\end{equation*}
where $\frac{1}{l}=\frac{1}{2}-\frac{1}{p}-\frac{1}{k}.$
By taking $R_2=J R_1 J^{-1}$ and $P=R_1 R_2,$ we have

\begin{equation}
\left\langle J, R_1: J^3=R_1^p=(R_1J)^{2k}=(R_1J^2)^{2l}=R_1 J R_1J^2 R_1 J R_1^{p-1} J^2 R_1^{p-1} J R_1^{p-1} J^2\right\rangle.
\end{equation}

A 3-fold Deligne-Mostow lattice of type three has the following presentation:

\begin{equation*}
\left\langle J, P, R_1, R_2: J^3=R_1^p=R_2^p=P^{3d}=(P^{-1}J)^k=Id, R_2=PR_1P^{-1} =J R_1 J^{-1}, P=R_1 R_2\right\rangle,
\end{equation*}
where $\frac{1}{d}=\frac{1}{2}-\frac{3}{p}$. By taking $R_2=J R_1 J^{-1}$ and $P=R_1 R_2,$ we obtain\begin{equation}
\left\langle J, R_1: J^3=R_1^p=(R_1J)^{2k}=(R_1 J R_1J^2)^{3d}=R_1 J R_1J^2 R_1 J R_1^{p-1} J^2 R_1^{p-1} J R_1^{p-1} J^2\right\rangle.
\end{equation}

A 3-fold Deligne-Mostow lattice of type four has the following presentation:

\begin{equation*}
\left\langle J, P, R_1, R_2: J^3=R_1^p=R_2^p=(P^{-1}J)^k=(R_2 R_1 J)^l=P^{3d}=Id, R_2=PR_1P^{-1} =J R_1 J^{-1}, P=R_1 R_2\right\rangle.
\end{equation*}
where $l$ and $d$ are defined as before. By taking $R_2=J R_1 J^{-1}$ and $P=R_1 R_2,$ we obtain
\begin{equation}
\left\langle J, R_1: J^3=R_1^p=(R_1J)^{2k}=(R_1J^2)^{2l}=(R_1JR_1J^2)^{3d}=R_1 J R_1J^2 R_1 J R_1^{p-1} J^2 R_1^{p-1} J R_1^{p-1} J^2\right\rangle.
\end{equation}

\section{Local rigidity of type preserving representations}

Given a presentation of a Deligne-Mostow lattice as above, we say that a representation of the lattice preserves generators types if the image of a complex reflection generator (respectively regular elliptic) is a complex reflection (respectively regular elliptic).  Note that in this definition the image of a complex reflection could be the identity and its order could be a factor of the order of the Deligne-Mostow generator.  That is, conjugacy classes might not be preserved.

\begin{theo}  Representations of Deligne-Mostow lattices with 3-fold symmetry into $\pgl(3,\C)$ preserving generator types in the above presentations are locally rigid and are classified in the following tables.\end{theo}

In the type one case, treated in \cite{ELA}, we know that local rigidity (with no conditions on types) holds for all groups with the exception of the three groups (4,4), (4,3) and (6,2).
We do not know whether rigidity holds for the groups of type two, three and four, when the generators are not of the prescribed type. 

In the following tables, we gather the information for each  3-fold  Deligne-Mostow lattice whose generators are of the same type as the original one, in the sense described above. 
The representations are defined with coefficients in a subfield of a cyclotomic extension of $\mathbb{Q},$ which depends on $(p,k)$ and for which we give the corresponding minimal polynomials. 
One of these images is the original Deligne-Mostow lattice and is contained in $PU(2,1)$. We write in boldface the field extension corresponding to this complex hyperbolic representation. 
We also compute the orbits of the Galois group action in the space of representations. 

Following Section 5 of \cite{ELA}, we will have representations coming from degenerate or non-degenerate configurations, according to the mutual positions of the fixed points sets of the generators. 
While irreducible representations can only come from non-degenerate configurations, reducible ones can come from either degenerate or non-degenerate configurations. In the Reducible column, the couple $(m,n)$ means that $m$ reducible representations come from non-degenerate configurations and $n$ come from degenerate ones. It is important to mention that, for each irreducible representation, there exists a unique (up to scalar multiplication) invariant Hermitian form whose signature might be different for different elements of the Galois orbit. The column Factors (which appears only in Table \ref{Tb:NonCompactDeligneMostowT2}) corresponds to the number of representations that contain in their kernel an element of the centralizer of the cusp group. One may interpret those representations as representations that factor through the orbifold fundamental group of a compactification of the complex hyperbolic orbifold defined by the Deligne-Mostow lattice. 

In Tables \ref{Tb:CompactT2}, \ref{Tb:DeligneMostowT3} and \ref{Tb:DeligneMostowT4} are the information of compact 3-fold of type two, three and four Deligne-Mostow lattices. Contrary to the compact 3-fold type one lattices, for those types there exist reducible representations coming from the non-degenerate configuration. 

\begin{table}[H]
\centering
\begin{tabular}{|c|c|p{1cm}|c|c|c|c|}
\hline
$(p,k)$ & Total & Galois Orbits &$\Q-$extension & Irreducible &Reducible & Factors\\ \hline
\multirow{5}{3em}{(6,4)} & \multirow{5}{3em}{46} & \multirow{5}{3em}{9} & $\mathbf{\Q(x^{12}-x^6+1)}$ & 2 & 0 & 1 \\ \cline{4-7} 
& & & $\Q(x^6-x^3+1)$ & 2 & 0 & 2 \\ \cline{4-7}
& & & $\Q(x^4-x^2+1)$ & 1 & 0  & 1 \\ \cline{4-7}
& & & $\Q(x^2-x+1)$ & 2 & 0  & 2 \\ \cline{4-7}
& & & $\Q$ & 1 & (0,1)  & 2  \\ \hline

\multirow{3}{3em}{(6,6)} & \multirow{3}{3em}{81} & \multirow{3}{3em}{30} & $\Q(x^6-x^3+1)$ & 6 & 0  & 4 \\ \cline{4-7}
& & & $\mathbf{\Q(x^2-x+1)}$ & 4 & (9,8)  & 15 \\ \cline {4-7}
& & & $\Q$ & 1 & (0,2)  & 3 \\ \hline

\end{tabular}
\caption{Representations of non-compact 3-fold Deligne-Mostow Lattices of Type two}
\label{Tb:NonCompactDeligneMostowT2}
\end{table}

\begin{table}[H]
\centering
%\begin{tabular}{|c|c|p{1cm}|c|c|c|}
\begin{tabular}{|c|c|c|c|c|c|}
\hline
$(p,k)$ & Total & Galois Orbits &$\Q-$extension & Irreducible & Reducible  \\ \hline
(3,7) & 36 & 1 & $\mathbf{\Q(x^{36}-x^{33}+x^{27}-x^{24}+x^{18}-x^{12}+x^9-x^3+1)}$ & 1 & 0  \\ \hline

\multirow{3}{3em}{(3,8)} & \multirow{3}{3em}{42} & \multirow{3}{3em}{3} & $\mathbf{\Q(x^{24}-x^{12}+1)}$ & 1 & 0  \\  \cline{4-6}
& & & $\Q(x^{12}-x^6+1)$ & 1 & 0 \\ \cline{4-6}
& & & $\Q(x^6-x^3+1)$ & 1 & 0 \\ \hline

\multirow{4}{3em}{(3,9)} & \multirow{4}{3em}{53} & \multirow{4}{3em}{16} & $\mathbf{\Q(x^{18}-x^6+1)}$ & 4 & 0  \\ \cline{4-6}
& & & $\Q(x^6-x^3+1)$ & 2 & 0 \\ \cline{4-6}
& & & $\Q(x^2-x+1)$ & 2 & (3,4) \\ \cline{4-6}
& & & $\Q$ & 0 & (0,1) \\\hline

\multirow{2}{3em}{(3,10)} & \multirow{2}{3em}{30} & \multirow{2}{3em}{2} & $\mathbf{\Q(x^{24}-x^{21}+x^{15}-x^{12}+x^9-x^3+1)}$ & 1 & 0  \\  \cline{4-6}
& & & $\Q(x^6-x^3+1)$ & 1 & 0 \\ \hline

\multirow{5}{3em}{(3,12)} & \multirow{5}{3em}{69}& \multirow{5}{3em}{18} & $\mathbf{\Q(x^{12}-x^6+1)}$ & 2 & 0  \\ \cline{4-6}
& & & $\Q(x^6-x^3+1)$ & 2 & 0 \\ \cline{4-6}
& & & $\Q(x^4-x^2+1)$ & 3 & 0  \\ \cline{4-6}
& & & $\Q(x^2-x+1)$ & 0 & (6,4)  \\ \cline{4-6}
& & & $\Q$ & 0 & (0,1) \\ \hline 

\multirow{4}{3em}{(4,5)} & \multirow{4}{3em}{37} & \multirow{4}{3em}{5} & $\mathbf{\Q(x^{16}+x^{14}-x^{10}-x^8-x^6+x^2+1)}$ & 2 & 0 \\ \cline{4-6}
& & & $\Q(x^8-x^7+x^5-x^4+x^3-x+1)$ & 1 & 0 \\ \cline{4-6}
& & & $\Q(x^4-x^3+x^2-x+1)$ & 1 & 0\\ \cline{4-6}
& & & $\Q$ & 0 & (0,1) \\ \hline 

\multirow{6}{3em}{(4,6)} & \multirow{6}{3em}{51} & \multirow{6}{3em}{12} & $\mathbf{\Q(x^{12}-x^6+1)}$ & 2 & 0  \\ \cline{4-6}
& & & $\Q(x^6-x^3+1)$ & 2 & 0  \\ \cline{4-6}
& & & $\Q(x^4-x^2+1)$ & 0 & (1,2)  \\ \cline{4-6}
& & & $\Q(x^2-x+1)$ & 1 & 0  \\ \cline{4-6}
& & & $\Q(x^2+1)$ & 0  & (1,0)  \\ \cline{4-6}
& & & $\Q$ & 1 & (0,2)  \\ \hline

\multirow{5}{3em}{(4,8)} & \multirow{5}{3em}{66} & \multirow{5}{3em}{17} & $\mathbf{\Q(x^8-x^4+1)}$ & 3 & 0  \\ \cline{4-6}
& & & $\Q(x^4-x^2+1)$ & 1 & (2,2)   \\ \cline{4-6}
& & & $\Q(x^4+1)$ & 3 & 0   \\ \cline{4-6}
& & & $\Q(x^2-x+1)$ & 2 & (2,0)  \\ \cline{4-6}
& & & $\Q$ & 1 & (0,1) \\ \hline

\multirow{3}{3em}{(5,4)} & \multirow{3}{3em}{36} & \multirow{3}{3em}{4} & $\mathbf{\Q(x^{16}+x^{14}-x^{10}-x^8-x^6+x^2+1)}$ & 1 & 0  \\ \cline{4-6}
& & & $\Q(x^8-x^7+x^5-x^4+x^3-x+1)$ & 2 & 0  \\ \cline{4-6}
& & & $\Q(x^4-x^3+x^2-x+1)$ & 1 & 0  \\  \cline{4-6}
\hline

\multirow{3}{3em}{(5,5)} & \multirow{3}{3em}{56} &\multirow{3}{3em}{10} & $\Q(x^8-x^7+x^5-x^4+x^3-x+1)$ & 3 & (1,2)  \\ \cline{4-6}
& & & $\mathbf{\Q(x^4-x^3+x^2-x+1)}$ & 3 & (1,0)  \\ \hline
 \end{tabular}
\caption{Representations of compact 3-fold Deligne-Mostow lattices of type two}
\label{Tb:CompactT2}
\end{table}

\begin{table}[H]
\centering
\begin{tabular}{|c|c|c|c|c|c|}
\hline
$(p,k)$ & Total & Galois Orbits &$\Q-$extension & Irreducible &Reducible  \\ \hline

\multirow{2}{3em}{$(7,2)$} & \multirow{2}{3em}{18} & \multirow{2}{3em}{2} & $\mathbf{\Q(x^{12}-x^{11}+x^9-x^8+x^6-x^4+x^3-x+1)}$ & 1 & 0  \\ \cline{4-6} 
& & & $\Q(x^6-x^5+x^4-x^3+x^2-x+1)$ & 1 & 0 \\ \hline

\multirow{6}{3em}{$(8,2)$} & \multirow{6}{3em}{26} & \multirow{6}{3em}{9} & $\mathbf{\Q(x^8-x^4+1)}$ & 1 & 0 \\ \cline{4-6}
& & & $\Q(x^4-x^2+1)$ & 0 & (1,0) \\ \cline{4-6}
& & & $\Q(x^4+1)$ & 1 & 0 \\ \cline{4-6}
& & & $\Q(x^2-x+1)$ & 1 & (0,2)  \\ \cline{4-6}
& & & $\Q(x^2+1)$ & 0 & (1,0) \\ \cline{4-6}
& & & $\Q$ & 1 & (0,1) \\ \hline

\multirow{2}{3em}{$(9,2)$} & \multirow{2}{3em}{24} & \multirow{2}{3em}{2} &  $\mathbf{\Q(x^{18}+x^9+1)}$ & 1 & 0 \\ \cline{4-6}
& & & $\Q(x^6-x^3+1)$ & 1 & 0 \\ \hline

\multirow{4}{3em}{$(10,2)$} & \multirow{4}{3em}{16} & \multirow{4}{3em}{5} & $\mathbf{\Q(x^8-x^7+x^5-x^4+x^3-x+1)}$ & 1 & 0  \\ \cline{4-6}
& & & $\Q(x^4-x^3+x^2-x+1)$ & 1 & 0 \\ \cline{4-6}
& & & $\Q(x^2-x+1)$ & 1 & 0 \\ \cline{4-6}
& & & $\Q$ & 1 & (0,1) \\ \hline

\multirow{5}{3em}{$(12,2)$} & \multirow{5}{3em}{37} & \multirow{5}{3em}{9} & $\mathbf{\Q(x^{12}-x^6+1)}$ & 1 & 0  \\ \cline{4-6}
& & & $\Q(x^6-x^3+1)$ & 1 & 0  \\ \cline{4-6}
& & & $\Q(x^4-x^2+1)$ & 0 & (1,2)  \\ \cline{4-6}
& & & $\Q(x^2-x+1)$ & 1 & 0  \\ \cline{4-6}
& & & $\Q(x^2+1)$ & 0 & (0,1) \\ \cline{4-6}
& & & $\Q$ & 1 & (0,1)  \\ \hline

\multirow{3}{3em}{$(18,2)$} & \multirow{3}{3em}{22} & \multirow{3}{3em}{4} & $\mathbf{\Q(x^{18}+x^9+1)}$ & 1 & 0  \\ \cline{4-6}
& & & $\Q(x^2-x+1)$ & 1 & 0 \\ \cline{4-6}
& & & $\Q$ & 1 & (0,1)  \\ \hline
\end{tabular}
\caption{Representations of 3-fold Deligne-Mostow Lattices of Type three}
\label{Tb:DeligneMostowT3}
\end{table}

\newpage 

\begin{table}[H]
\centering
\begin{tabular}{|c|c|c|c|c|c|}
\hline
$(p,k)$ & Total & Galois Orbits &$\Q-$extension & Irreducible &Reducible \\ \hline
$(7,3)$ & 36 &1&$\mathbf{\Q(x^{24}-x^{20}+x^{16}-x^{12}+x^8-x^4+1)}$ &1 & 0 \\ \hline

\multirow{3}{3em}{$(8,3)$} & \multirow{2}{3em}{46} & \multirow{2}{3em}{5} & $\mathbf{\Q(x^{24}-x^{12}+1)}$ & 1 & 0 \\ \cline{4-6}
& & & $\Q(x^{12}-x^6+1)$ & 1 & 0\\ \cline{4-6}
& & & $\Q(x^6-x^3+1)$ & 1 & 0\\ \cline{4-6}
& & & $\Q$ & 0 & (0,2)\\ \hline

\multirow{6}{3em}{$(8,4)$} & \multirow{6}{3em}{74} & \multirow{6}{3em}{19} & $\mathbf{\Q(x^8-x^4+1)}$ & 2 & (0,1) \\ \cline{4-6}  
& & & $\Q(x^4-x^2+1)$ & 1 & (2, 2)\\ \cline{4-6} 
& & & $\Q(x^4+1)$ & 2 & (1,2)\\ \cline{4-6}
& & & $\Q(x^2-x+1)$ & 1 & 0 \\ \cline{4-6}
& & & $\Q(x^2+1)$ & 1 & (2,0)\\ \cline{4-6}
& & & $\Q$ &1 & (0,1)\\ \hline 

\multirow{4}{3em}{$(9,3)$} & \multirow{4}{3em}{53} & \multirow{4}{3em}{11} & $\mathbf{\Q(x^{18}+x^9+1)}$ & 2 & 0\\ \cline{4-6} 
& & & $\Q(x^6-x^3+1)$ & 1 & 0\\ \cline{4-6}
& & & $\Q(x^2-x+1)$ & 0 & (3,4)\\ \cline{4-6}
& & & $\Q$ & 0 & (0,1)\\  \hline

\multirow{3}{3em}{$(10,3)$} & \multirow{3}{3em}{32} & \multirow{3}{3em}{4} & $\mathbf{\Q(x^{24}+x^{21}-x^{15}-x^{12}-x^9+x^3+1)}$ & 1 & 0\\ \cline{4-6}
& & & $\Q(x^6-x^3+1)$ & 1 & 0 \\ \cline{4-6}
& & & $\Q(x^2+x+1)$ & 0 & (0,2) \\ \hline 

\multirow{3}{3em}{$(10,5)$} & \multirow{3}{3em}{65} & \multirow{3}{3em}{11} & $\Q(x^8-x^7+x^5-x^4+x^3-x+1)$ & 3 & (2,2)\\ \cline{4-6}
& & & $\mathbf{\Q(x^4-x^3+x^2-x+1)}$ & 1 & (2,0)\\ \cline{4-6}
& & & $\Q$ & 0 & (0,1)\\ \hline

\multirow{5}{3em}{$(12,3)$} & \multirow{5}{3em}{61} & \multirow{5}{3em}{14} & $\Q(x^{12}-x^{6}+1)$ & 2 & 0  \\ \cline{4-6}
& & & $\Q(x^6-x^3+1)$ & 2 & 0\\ \cline{4-6}
&  & & $\mathbf{\Q(x^4-x^2+1)}$ & 3 & 0\\ \cline{4-6}
& & & $\Q(x^2-x+1)$ &0 & (6,0) \\ \cline{4-6}
& & & $\Q$ & 0 & (0,1) \\ \hline

\multirow{6}{3em}{$(12,6)$} & \multirow{6}{3em}{53} & \multirow{6}{3em}{15} & $\mathbf{\Q(x^{12}-x^6+1)}$ & 2 & 0\\ \cline{4-6}
& & & $\Q(x^6-x^3+1)$ & 1 & 0\\ \cline{4-6}
& & & $\Q(x^4-x^2+1)$ & 1 & 0\\ \cline{4-6}
& & & $\Q(x^2-x+1)$ & 1 & (0,6)\\ \cline{4-6}
& & & $\Q(x^2+1)$ & 1 & 0\\ \cline{4-6}
& & & $\Q$ & 1 & (0,2)\\ \hline

\multirow{4}{3em}{$(18,3)$} & \multirow{4}{3em}{58} & \multirow{4}{3em}{14} & $\mathbf{\Q(x^{18}-x^9+1)}$ & 2 & 0 \\ \cline{4-6} 
& & & $\Q(x^6-x^3+1)$ & 1 & 0\\ \cline{4-6} 
& & & $\Q(x^2-x+1)$ & 0 & (3,6)\\ \cline{4-6} 
& & &  $\Q$ & 0 & (0,2)\\ \hline
\end{tabular}
\caption{Representations of 3-fold Deligne-Mostow Lattices of Type four}
\label{Tb:DeligneMostowT4}
\end{table}

\addcontentsline{toc}{section}{\refname}
\bibliographystyle{alpha}
\bibliography{biblio}
%\nocite{*}

\begin{flushleft}
  \textsc{E. Falbel\\
  Sorbonne Universit\'e and Universit\'e de Paris, CNRS, INRIA, IMJ-PRG, \\ F-75005 Paris, France
  \\}
 \verb|elisha.falbel@imj-prg.fr|
 \end{flushleft}
 
  \begin{flushleft}
  \textsc{I. Pasquinelli\\
  University of Bristol  \\
School of Mathematics\\
Fry Building 
Woodland Road 
Bristol BS8 1UG 
UK
 \\}
 \verb|irene.pasquinelli@bristol.ac.uk|
 \end{flushleft}
 
 \begin{flushleft}
  \textsc{A. Ucan-Puc\\
  Sorbonne Universit\'e and Universit\'e de Paris, CNRS, INRIA, IMJ-PRG,\\
   F-75005 Paris, France
\\}
 \verb|alejandro.ucan-puc@imj-prg.fr|
 \end{flushleft}

\end{document}